\documentclass[12pt,reqno]{article}
\usepackage{amssymb}
\usepackage{amsmath}
\usepackage{amsthm}
\usepackage{amssymb}
\usepackage{amsmath}
\usepackage{amsthm}
\newcommand{\br}[3]{{$#1$}$\lower4pt\hbox{$\tp\atop\raise4pt \hbox{$\scriptscriptstyle{#2}$}$} ${$#3$}}
\newcommand{\tw}[3]{{$#1$}${\,\scriptscriptstyle {#2}}\atop\raise9pt\hbox{$\scriptstyle\tp$} ${$#3$}}
\newcommand{\ttps}[2]{{#1}\raise5pt\hbox{$\lower12pt\hbox{$\scriptstyle\tp$}\atop \lower0pt\hbox{$\tilde\;$}$}\raise4.5pt\hbox{${\scriptstyle{#2}}$}}
\newcommand{\otp}{\tilde \tp}
\newcommand{\otps}{\tilde \tp}
\newcommand{\ttp}{{\lower12pt\hbox{$\tp$}\atop \hbox{$\tilde\;$}}}

\newcommand{\btr}{\raise1.2pt\hbox{$\scriptstyle\blacktriangleright$}\hspace{2pt}}
\newcommand{\id}{\mathrm{id}}

\newcommand{\tp}{\otimes}

\newcommand{\A}{\mathcal{A}}
\newcommand{\B}{\mathcal{B}}

\newcommand{\N}{\mathbb{N}}
\newcommand{\C}{\mathbb{C}}
\newcommand{\R}{\mathbb{R}}
\newcommand{\Z}{\mathbb{Z}}

\newcommand{\Ha}{\mathcal{H}}
\newcommand{\Ru}{\mathcal{R}}
\newcommand{\Q}{\mathcal{Q}}
\newcommand{\U}{\mathcal{U}}
\newcommand{\F}{\mathcal{F}}
\newcommand{\ve}{\varepsilon}
\newcommand{\gm}{\gamma}
\newcommand{\End}{\mathrm{End}}
\newcommand{\tr}{\triangleright}
\newcommand{\tl}{\triangleleft}

\newcommand{\g}{\mathfrak{g}}
\newcommand{\K}{\mathcal{K}}
\newcommand{\Ta}{\mathcal{T}}
\newcommand{\KK}{\mathbf{k}}

\newcommand{\I}{\mathrm{I}}
\newcommand{\si}{\sigma}
\newcommand{\be}{\begin{eqnarray}}
\newcommand{\ee}{\end{eqnarray}}
\newtheorem{thm}{Theorem}
\newtheorem{propn}[thm]{Proposition}

\newtheorem{corollary}[thm]{Corollary}

\theoremstyle{definition}
\newtheorem{parag}{}[section]
\newtheorem{remark}[thm]{Remark}
\newtheorem{definition}{Definition}
\newtheorem{example}{Example}

\begin{document}
\title{On universal solution to reflection equation\footnote{This research is partially supported
by the Israel Academy of Sciences grant no. 8007/99-01 and by the RFBR grant no. 02-01-00085.}}
\author{J. Donin$^\dag$,
P. P. Kulish$^{\ddag}$, and A. I.  Mudrov$^\dag$}
%
\setcounter{page}{1}
\maketitle
\begin{center}
{\small \dag \hspace{5pt}Department of Mathematics, Bar Ilan University,
52900 Ramat Gan, Israel.}
\\
{\small \ddag \hspace{5pt}St.-Petersburg Department of Steklov Mathematical Institute,
Fontanka 27, \\ \indent\hspace{-7.5cm} {191011 St.-Petersburg, Russia}}
\end{center}
\begin{abstract}
For a given  quasitriangular Hopf algebra $\Ha$ we study relations between
the braided group $\tilde \Ha^*$ and Drinfeld's twist.
We show that the braided bialgebra structure of $\tilde \Ha^*$ is naturally described by
means of twisted tensor powers of $\Ha$ and their module algebras.
We introduce universal solution to the reflection equation (RE) and deduce a fusion prescription
for  RE-matrices.\\[16pt]
{\small \underline{Key words}: reflection equation, twist, fusion procedure.\\
\underline{AMS classification codes}: 17B37, 16W30.}
\end{abstract}
\section{Introduction}
There are two important algebras in the quantum group theory, the Faddeev-Reshetikhin-Takhtajan (FRT)
and reflection equation (RE) algebras. As was shown in \cite{DM}, they are related by a twist of
the underlying quasitriangular Hopf algebra $\Ha$ squared. This twist transforms any bimodule over
$\Ha$ to an $\Ha^{\otps2}$-module, so there is an analog of the dual algebra
$\Ha^*$, an $\Ha^{\otps2}$-module algebra $\tilde \Ha^*$ . This algebra turns out to be isomorphic to the Hopf algebra in the quasitensor category
of $\Ha$-modules, \cite{Mj}, thus equipped with an additional structure of braided coalgebra.
In the present paper we develop further the approach of \cite{DM} studying this structure from
the twist point of view.
The algebra $\tilde \Ha^*$ is a module over the twisted tensor square $\Ha^{\otps2}$. This twist
can be extended to all tensor powers of $\Ha$ giving $\Ha^{\otps n}$, $n=0,1,\ldots,$ where
the zero power is set to be the ring of scalars with the natural Hopf algebra structure.
It appears  natural to consider modules over all Hopf algebras $\Ha^{\otps n}$ simultaneously. They form
a monoidal category with respect to twisted tensor product, and the iterated comultiplications
$\Delta^n\colon\Ha\to \Ha^{\otps n}$ induce a functor from that category
to the quasitensor category of $\Ha$-modules.
The iterated braided coproducts $\Delta^n_{\tilde\Ha^*}$ on $\tilde\Ha^*$ take their values in module algebras over $\Ha^{\otps 2n}$, $n=2,3,\ldots$

We introduce universal K-matrix, $\K\in\Ha\tp \tilde\Ha^*$, satisfying the
characteristic equation
$$(\Delta\tp \id)(\K)=\Ru^{-1}\K_1\Ru\K_2.$$
This equation implies the
"abstract" reflection equation
$$\Ru_{21} \K_1 \Ru \K_2 =  \K_2 \Ru_{21} \K_1 \Ru$$
in $\Ha\tp \Ha\tp \tilde\Ha^*$ involving the universal R-matrix.
We prove that $\K$ is equal to the canonical element of $\Ha\tp \tilde\Ha^*$
with respect to the Hopf pairing between $\Ha$ and $\Ha^*$, which coincides with $\tilde\Ha^*$
as a linear space.
The reason for considering this object is the same as for the universal
R-matrix: it gives a solution to the matrix RE in  every representation of
$\Ha$.

A fusion procedure for R-matrices was
the first contact of the Yang-Baxter
equation with algebraic structures [KRS].
Later it was related with properties
of the universal R-matrix in the theory
of quantum groups [Dr1]. Although an analogous
fusion procedure was considered for matrix
solutions to the reflection equation [MeN, KSkl],
to our knowledge no universal element was
proposed  in any algebraic approach to the RE
(see, e.g., [DeN] and references therein).
Using the characteristic equation on the universal element $\K$, we formulate a version of
fusion procedure for RE matrices, suggesting an algorithm of tensoring
RE matrices.

The paper is organized as follows. Section \ref{sQHA} contains a summary on
Hopf algebras, universal R-matrix, and twist. There we introduce RE matrices
and prove an auxiliary proposition about their restriction to submodules.
In Section \ref{TTP} we consider a special type of twist
of tensor product Hopf algebras and their modules when the twisting cocycle satisfies the bicharacter identities.
We show that the braided tensor product of module algebras is a twisted tensor product. In Section \ref{sTTpwr} we extend this construction
to higher tensor powers of  Hopf algebras. In Section \ref{sREdual} we focus
our study on the RE dual $\tilde \Ha^*$ to a quasitriangular Hopf algebra $\Ha$. There, we introduce
the universal RE matrix  $\K$ and deduce the characteristic equation for it.
Using the bicharacter twist and the universal K-matrix, we study the braided bialgebra structure
$\tilde \Ha^*$ in Section \ref{sBCT}.
Section \ref{sFusion} is devoted to the central result of the paper,
a fusion procedure for RE matrices. Appendix contains the proof of Theorem \ref{fusion_thm}.


\section{Preliminaries}
\label{sQHA}
\begin{parag}
Let $\KK$ be a commutative algebra over a field
of zero characteristic\footnote{One may think of $\KK$ as either $\C$ or $\C[[h]]$. In the latter case
all algebras are assumed complete in the $h$-adic topology.}. Let $\Ha$ be a quasitriangular Hopf  algebra\footnote{
For a guide in quasitriangular Hopf algebras, the reader is
referred to original Drinfeld's report, \cite{Dr1}, and to one of
the textbooks, e.g. \cite{ChPr} or \cite{Mj}.}
over $\KK$,
with the coproduct $\Delta$, counit $\ve$, antipode $\gm$, and
the universal R-matrix $\Ru\in \Ha\tp \Ha$.
By definition, it satisfies following Drinfeld's equations, \cite{Dr1}:
\be{}
(\Delta \tp \id)(\Ru)=\Ru_{13}\Ru_{23}, \quad (\id \tp\Delta)(\Ru)=\Ru_{13}\Ru_{12},
\label{bch}
\ee{}
and, for any $x\in \Ha$,
\be{}
\quad\Ru\Delta(x)=\Delta^{op}(x)\Ru.
\label{flip}
\ee{}
We adopt the usual convention of marking tensor components
with subscripts indicating the supporting tensor factors.
The subscript $op$ will
stand for the opposite multiplication while the superscript $op$
for the opposite coproduct. The standard symbolic notation
with suppressed summation
is used for the coproduct, $\Delta(x)=x_{(1)}\tp x_{(2)}$.

By $\N$ we assume a set of non-negative integers, i.e., including zero.
By $\Delta^n$, $n\in \N$,  we denote the n-fold coproduct, $\Delta^n\colon \Ha\to \Ha^{\tp n}$, setting
\be
\label{iterated}
\Delta^0=\ve,\quad \Delta^1=\id,\quad \Delta^2=\Delta,\quad \Delta^3=(\Delta\tp \id)\circ \Delta,
\quad \ldots
\ee
Here and further on we view  the ring of scalars $\KK$ as equipped with the structure
of Hopf algebra over $\KK$.
It is then convenient to put $\Ha^{\tp n}=\KK$ for $n=0$.

The  antipode $\gm$ is treated as a Hopf algebra  isomorphism
between $\Ha_{op}$ and $\Ha^{op}$. The universal R-matrix defines two
homomorphisms from the dual coopposite Hopf algebra $\Ha^{*op}$ to $\Ha$:
\be{}
\Ru^\pm(a) = \langle a,\Ru^\pm_{1}\rangle \Ru^\pm_{2},
\quad a\in \Ha^{*op},\quad \mbox{where}\quad
\Ru^+ = \Ru\quad \mbox{and }\quad\Ru^- = \Ru^{-1}_{21}.
\label{R+-}
\ee{}
\end{parag}
\begin{parag}
Twist of a Hopf algebra $\Ha$ by a cocycle $\F$ is a Hopf algebra
$\tilde \Ha$ whose comultiplication  is obtained from $\Delta$ by a similarity transformation, \cite{Dr2},
$$\tilde\Delta(x)=\F^{-1}\Delta(x)\F,\quad x\in \tilde\Ha.$$
To preserve coassociativity, it is sufficient for the element $\F\in \Ha\tp \Ha$ to satisfy the cocycle
equation
\be
\label{cc}
(\Delta \tp \id)(\F)\F_{12}=(\id\tp\Delta)(\F)\F_{23}
\ee
and the normalization condition
\be
\label{norm}
(\ve\tp \id)(\F)=( \id\tp\ve)(\F)=1\tp1.
\ee
The multiplication and counit in $\tilde\Ha$ remain the same as in $\Ha$.
The antipode changes
by a similarity transformation,
$\tilde \gm(x) = u^{-1}\gm(x)u$ with $u=\gm(\F_1)\F_2 \in \Ha$,
and the universal R-matrix is
\be{}
\label{Rtwisted}
\tilde \Ru = \F^{-1}_{21}\Ru \F.
\ee{}
Twist establishes an equivalence relation among Hopf algebras, so
we call $\tilde \Ha$ and $\Ha$ twist-equivalent. We use notation
$\tilde \Ha\stackrel{\F}{\sim}\Ha$; then $\Ha\stackrel{\F^{-1}}{\sim}\tilde \Ha$.

Recall that a left $\Ha$-module algebra $\A$ is an associative
algebra over $\KK$ endowed with the left action $\btr$ of $\Ha$ such that
the multiplication map $\A\tp\A\to \A$
is $\Ha$-equivariant.
Given two $\Ha$-module algebras $\A$ and $\B$, there exists
their braided tensor product,
 \cite{Mj}. This
is an $\Ha$-module algebra coinciding with $\A\tp \B$ as a linear
space and endowed with multiplication
\be{}
\label{br_t_pr}
\bigl(a_1\tp b_1\bigr)\bigl(a_2\tp b_2\bigr)=a_1 (\Ru_2 \btr a_2) \tp
\bigl(\Ru_1 \btr b_1\bigr)\:b_2,
\ee{}
for $a_i\in \A$ and $b_i\in \B$, $i=1,2$.

An $\Ha$-module algebra $\A$ becomes a left $\tilde\Ha$-module algebra $\tilde \A$, $\tilde \A\stackrel{\F}{\sim}\A$,
when equipped with the multiplication
\be
a\circ b = (\F_1\btr a)(\F_2\btr b), \quad a,b \in \A.
\ee
The action of $\tilde\Ha$ on $\tilde \A$ is that of  $\Ha$ on $\A$,
having in mind $\tilde\Ha \simeq\Ha$ as associative algebras.
\end{parag}

\begin{parag}[Reflection equation]
Let $(\rho,V)$ be a representation of $\Ha$ on a module $V$, 
and
$R$ the image of the universal R-matrix, $R=(\rho\tp \rho)(\Ru)\in \End^{\tp 2}(V)$.
Let $\A$ be an associative algebra. An element $K\in \End(V)\tp \A$ is said to be
a solution to (constant) RE or an (constant) RE matrix in representation $\rho$ with coefficients
in $\A$ if
\be
\label{re}
R_{21}K_1RK_2=K_2R_{21}K_1R.
\ee
This equation is supported in $\End^{\tp 2}(V)\tp \A$, and the subscripts
label the components belonging to the different tensor factors $\End(V)$.
Let us prove the following elementary proposition.
\begin{propn}
\label{restriction}
Let $K\in \End(V)\tp \A$  be an RE matrix and suppose the $\Ha$-module $V$ is semisimple.
Let $\rho_0$ be a representation of $\Ha$ on the submodule $V_0\subset V$
and $V_0\stackrel{\iota}{\to} V\stackrel{\pi}{\to} V_0$ the intertwiners.
Then, the matrix $K_0=\pi K\iota \in \End(V_0)\tp \A$ is a solution to the RE in the representation $\rho_0$.
\end{propn}
\begin{proof}
Multiply equation (\ref{re})
by $\pi\tp\pi$ from the left and by $\iota\tp\iota$ from the right and use
the intertwining formulas $\pi \rho(x)=\rho_0(x)\pi$, $\rho(x)\iota=\iota\rho_0(x)$
valid for every $x\in \Ha$.
The result will be the RE on the matrix $K_0$  with
the R-matrix $(\rho_0\tp\rho_0)(\Ru)$.
\end{proof}

\end{parag}
\section{Twisted tensor product of Hopf algebras}
\label{TTP}
\begin{parag}
In this section we recall a  specific case of twist as applied to tensor products of
Hopf algebras, \cite{RS}. This construction is of particular importance for our consideration.
Denote by  $\Ha'$ the tensor product $\Ha^{\{1\}}\tp \Ha^{\{2\}}$ of two Hopf algebras.
It is equipped with the standard tensor product multiplication and comultiplication.
An element $\F\in \Ha^{\{2\}}\tp \Ha^{\{1\}}$ may be viewed as that from
$\Ha'\tp\Ha'$  via the embedding
$$\F\in\bigl(1\tp\Ha^{\{2\}}\bigr)\tp \bigl(\Ha^{\{1\}}\tp 1\bigr)\subset\Ha'\tp\Ha'.$$
If $\F$ satisfies the bicharacter identities
\begin{equation}
\begin{array}{ccc}
(\Delta^{\{2\}}\tp\id)(\F)&=&\F_{13}\F_{23}\in \Ha^{\{2\}}\tp\Ha^{\{2\}}\tp\Ha^{\{1\}}
,\\[6pt]
(\id\tp\Delta^{\{1\}})(\F)&=&\F_{13}\F_{12}\in \Ha^{\{2\}}\tp\Ha^{\{1\}}\tp\Ha^{\{1\}},
\end{array}
\label{ttp}
\end{equation}
 then it fulfills cocycle condition (\ref{cc}) in $\Ha'^{\tp 2}$ and condition (\ref{norm}).
\begin{definition}
{\em Twisted tensor product} \tw{\Ha^{\{1\}}}{\F}{\Ha^{\{2\}}}  of two Hopf
algebras $\Ha^{\{i\}}$, $i=1,2$, is the twist of $\Ha^{\{1\}}\tp \Ha^{\{2\}}$ with a bicharacter
cocycle $\F$ satisfying
(\ref{ttp}).
\end{definition}
\noindent
\begin{propn}
\label{projections}
The maps
$$\id \tp \ve^{\{2\}}\colon \mbox{\tw{\Ha^{\{1\}}}{\F}{\Ha^{\{2\}}}}\to \Ha^{\{1\}}
,\quad
  \ve^{\{1\}} \tp \id\colon \mbox{\tw{\Ha^{\{1\}}}{\F}{\Ha^{\{2\}}}}\to \Ha^{\{2\}}
$$
are Hopf algebra homomorphisms.
\end{propn}
\begin{proof}
This is a corollary of the equalities $(\ve^{\{2\}}\tp \id)(\F)=1\tp 1=(\id \tp\ve^{\{1\}})(\F)$
which follow from (\ref{ttp}).
\end{proof}
\end{parag}
\begin{parag}
\label{ttsH}
A particular case of twisted tensor product is the twisted tensor square \tw{\Ha}{\Ru}{\Ha}
of a quasitriangular Hopf algebra $\Ha$.
The universal R-matrix satisfies bicharacter conditions (\ref{ttp}) by virtue of
(\ref{bch}).
If one takes the universal R-matrix of $\Ha\tp\Ha $ equal to $\Ru^-_{13}\Ru^+_{24}$, then
formula (\ref{Rtwisted}) gives the universal R-matrix of \tw{\Ha}{\Ru}{\Ha}:
\be{}
\label{ttpRm}
\Ru^{-}_{14}\Ru^-_{13}\Ru^+_{24} \Ru_{23}=
\Ru^{-1}_{41}\Ru^{-1}_{31}\Ru _{24} \Ru_{23}
\in
(\mbox{\tw{\Ha}{\Ru}{\Ha}})\tp(\mbox{\tw{\Ha}{\Ru}{\Ha}}).
\ee{}
\begin{propn}
\label{main}
The Hopf algebra \tw{\Ha}{\Ru}{\Ha} is twist-equivalent to
$\Ha^{op}\tp\Ha$ via the cocycle $\F=\Ru_{13}\Ru_{23}\in (\Ha^{op}\tp\Ha)\tp(\Ha^{op}\tp\Ha)$.
\end{propn}
\begin{proof}
This twist is a composition of the two,
$\Ha\tp \Ha \stackrel{\Ru_{13}}{\sim} \Ha^{op}\tp \Ha$ and
$\mbox{\tw{\Ha}{\Ru}{\Ha}} \stackrel{\Ru_{23}}{\sim} \Ha\tp \Ha$.
\end{proof}
\end{parag}

\begin{parag}[Twisted tensor product of module algebras]
\label{skew_tp}
Let $\Ha^{\{i\}}$, $i=1,2$, be  Hopf algebras and  $\A^{\{i\}}$ their left
module algebras.
Consider the twisted tensor product
$\Ha'=\mbox{\tw{\Ha^{\{1\}}}{\F}{\Ha^{\{2\}}}}$
by a bicharacter $\F$.
There is a left $\Ha'$-module algebra, \br{\A^{\{1\}}}{\F}{\A^{\{2\}}}, which is built on
$\A^{\{1\}}\tp \A^{\{2\}}$ as an $\Ha'$-module and contains
$\A^{\{1\}}\tp 1$ and  $1\tp \A^{\{2\}}$ as invariant subalgebras.
The multiplication is given by
\be{}
\label{br_pr}
\bigl(a^{\{1\}}_1\tp a^{\{2\}}_1\bigr)\bigl(a^{\{1\}}_2\tp a^{\{2\}}_2\bigr)=a^{\{1\}}_1 (\F_2 \btr a^{\{1\}}_2) \tp
\bigl(\F_1 \btr a^{\{2\}}_1\bigr)\:a^{\{2\}}_2,
\ee{}
where $a^{\{i\}}_j\in \A^{\{i\}}$.
This multiplication satisfies the permutation
relation
\be{}
\label{br_rel} (1\tp a^{\{2\}})(a^{\{1\}}\tp 1)=\bigl(\F_2 \btr a^{\{1\}} \tp 1\bigr)
\bigl(1\tp \F_1 \btr a^{\{2\}}\bigr), \quad a^{\{i\}}\in \A^{\{i\}}.
\ee{}
\begin{remark}
\label{btp}
A particular case of this construction, when $\Ha^{\{1\}}=\Ha^{\{2\}}=\Ha$ and $\F=\Ru$
(cf. Subsection \ref{ttsH}),
gives the braided tensor product  (\ref{br_t_pr}) of $\Ha$-module algebras, which acquires
an $\Ha$-module algebra structure because of the Hopf algebra homomorphism $\Delta\colon\Ha\to \mbox{\tw{\Ha}{\Ru}{\Ha}}$.
\end{remark}
\end{parag}
\newcommand{\otpb}{\mbox{$\raise1pt \hbox{
$\hspace{-6pt}\scriptscriptstyle\bigcirc \hspace{-8.9pt}\otimes\hspace{1.5pt}$}$}}

\section{Twisted tensor powers $\Ha^{\otp n}$}
\label{sTTpwr}
\begin{parag}
The purpose of this section is to extend the collection of $\Ha$-modules
by modules over certain Hopf algebras related to $\Ha$. The reason for that
is the observation that the braided tensor product of two $\Ha$-module algebras, is a particular
case of the twisted tensor product of modules (cf. Remark \ref{btp}) and
thus a module over \tw{\Ha}{\Ru}{\Ha}.
\end{parag}
\begin{parag}
Consider the set $\mathfrak{H}(\Ha)=\cup_{\Ha'}\mathrm{Hom}(\Ha,\Ha')$ of Hopf algebra homomorphisms,
where $\Ha'$ runs over Hopf algebras.
We treat the ring $\KK$ of scalars as a Hopf algebra over $\KK$, so $\ve\in\mathfrak{H}(\Ha)$.
Let $\phi^{\{i\}}\in \mathfrak{H}(\Ha)$, $i=1,2$,  be two homomorphisms
from $\Ha$ to Hopf algebras $\Ha^{\{i\}}$.
Put the bicharacter $\F\in \Ha^{\{2\}}\tp\Ha^{\{1\}}$ to be the image of the universal
R-matrix, $\F=\bigl(\phi^{\{2\}}\tp\phi^{\{1\}}\bigr)(\Ru)$.

\begin{propn}
The map
\be
\phi^{\{1\}}\otp \phi^{\{2\}}=\bigl(\phi^{\{1\}}\tp\phi^{\{2\}}\bigr)\circ \Delta\colon \Ha\to \mbox{\tw{\Ha^{\{1\}}}{\F}{\Ha^{\{2\}}}}
\ee
is a Hopf algebra homomorphism.
\end{propn}
\begin{proof}
This follows from the readily verified fact that the map $\Delta\colon \Ha\to \mbox{\tw{\Ha}{\Ru}{\Ha}}$ is a Hopf
algebra embedding.
\end{proof}
\noindent
\begin{propn}
The operation $\otp$ makes the set $\mathfrak{H}(\Ha)$ a monoid
with $\ve$ being the neutral element. The operation $\otp$ is natural:
if $f^{\{i\}}\colon\Ha^{\{i\}}\to \Ha^{\{i\}}_1$ are Hopf algebra homomorphisms, then
$$
\bigl(f^{\{1\}}\circ \phi^{\{1\}}\bigr)\otp \bigl(f^{\{2\}}\circ \phi^{\{2\}}\bigr)=
\bigl(f^{\{1\}}\tp f^{\{2\}}\bigr)\circ \bigl(\phi^{\{1\}}\otp \phi^{\{2\}}\bigr).
$$
\end{propn}
\begin{proof}
It follows immediately from coassociativity of $\Delta$ that the operation $\otp$ is associative.
It is obviously natural because so is the construction of twisted tensor product.
The counit map is the identity due to the arguments used in the proof of Proposition \ref{projections}.
\end{proof}
We will use the notation $\Ha^{\{1\}}\otp\Ha^{\{2\}}$ for \tw{\Ha^{\{1\}}}{\F}{\Ha^{\{2\}}} when
the specific form of homomorphisms $\phi^{\{i\}}$ is clear from the context.
In particular, we consider the set  $\{\id^{\otps n}|\; n\in \N\}$,  where $\id^0=\ve$, and
the corresponding twisted tensor powers $\{\Ha^{\otps n}|\; n\in \N\}$.
 The projections $\Ha^{\otps n}\to \Ha^{\otps k}$ obtained by applying $\ve$ to some tensor
factors $n-k$ times are Hopf algebra homomorphisms.
Also, for every $k\in\N$, the coproduct $\Delta^k$ applied to any component of $\Ha^{\otps n}$
induces a Hopf algebra homomorphism from $\Ha^{\otps n}$ to $\Ha^{\otps (n+k)}$
\begin{parag}
Let $\A^{\{i\}}$ be two left $\Ha$-module algebras.
By $\A^{\{1\}}\ttp\A^{\{2\}}$ we denote the $\Ha^{\otps2}$-module algebra
\br{\A^{\{1\}}}{\Ru}{\A^{\{2\}}}, see Subsection \ref{skew_tp}. Suppose  $\A^{\{i\}}$  are module algebras
over Hopf algebras $\Ha^{\{i\}}$ and the $\Ha$-module structures on $\A^{\{i\}}$ are
induced by homomorphisms $\phi^{\{i\}}\colon \Ha \to \Ha^{\{i\}}$.
Obviously, the algebra $\A^{\{1\}}\ttp \A^{\{2\}}$ is isomorphic
to \br{\A^{\{1\}}}{\F}{\A^{\{2\}}}, with $\F=\bigl(\phi^{\{2\}}\tp\phi^{\{1\}}\bigr)(\Ru)$,
and the $\Ha^{\otps 2}$-module structure on it is induced by
the homomorphism
$\phi^{\{1\}}\tp \phi^{\{2\}}\colon \Ha^{\otps2} \to \mbox{\tw{\Ha^{\{1\}}}{\F}{\Ha^{\{2\}}}}$
of Hopf algebras.
\begin{example}
The operation $\ttp$ is associative and allows to introduce twisted tensor powers of
an $\Ha$-module algebra $\A$. By induction, we can define the $n$-th tensor power $\ttps{\A}{n}$,
setting $\ttps{\A}{n}=\ttps{\A}{(n-1)}\ttp \A$.
\end{example}
\begin{example}
Let $\Ha'$ be another Hopf algebra such that there exists a Hopf algebra
homomorphism $\phi\colon \Ha \to \Ha'$. Suppose that the algebra $\A$ from the previous example is also an $\Ha'$-module algebra
and the $\Ha$-module structure on $\A$ is induced by $\phi$.
Then, $\ttps{\A}{n}$ is a module over $\Ha'^{\otps n}$,  $\Ha^{\otps n}$, and  $\Ha$ (but not over $\Ha'$
since, in general, there is no natural map from $\Ha'$ to  $\Ha'^{\otps n}$,
which is twisted over $\Ha$).
\end{example}
\end{parag}
\end{parag}
\section{RE dual $\tilde \Ha^*$}
\label{sREdual}
\begin{parag}[RE twist]
Consider a quasitriangular Hopf algebra $\Ha'$ with the universal R-matrix $\Ru'$ and a left
$\Ha'$-module algebra $\A$.
\begin{definition}
$\A$ is called {\em quasi-commutative} if for
any $a,b\in \A$
\be{}
\label{crel}
(\Ru'_2\btr b)(\Ru'_1\btr a) = a b.
\ee{}
\end{definition}
\noindent
The following statement is elementary.
\begin{propn}
\label{tw-com}
Let $\Ha'$ be a quasitriangular Hopf algebra and
$\tilde \Ha'\stackrel{\F}{\sim}\Ha'$.
If $\A$ is a quasi-commutative $\Ha'$-module algebra,
then the twisted algebra $\tilde \A$,  $\tilde\A\stackrel{\F}{\sim}\A$,
is a quasi-commutative $\tilde \Ha'$-module algebra.
\end{propn}
\begin{proof}
Clear.
\end{proof}
\begin{corollary}
\label{qctp}
Let $\A^{\{i\}}$ be quasi-commutative module algebras over
Hopf algebras $\Ha^{\{i\}}$, where $i=1,2$. Let $\F\in \Ha^{\{2\}}\tp \Ha^{\{1\}}$ be a bicharacter
and $\Ha'$ the twisted tensor product \tw{\Ha^{\{1\}}}{\F}{\Ha^{\{2\}}}.
Then, the twisted tensor product \br{\A^{\{1\}}}{\F}{\A^{\{2\}}} is
a quasi-commutative $\Ha'$-module algebra.
\end{corollary}
\begin{proof}
The ordinary tensor product $\A^{\{1\}}\tp \A^{\{2\}}$
is a quasi-commutative $\Ha^{\{1\}}\tp \Ha^{\{2\}}$-module algebra.
Now apply Proposition \ref{tw-com}.
\end{proof}
\begin{example}
The Hopf dual $\Ha^*$ is an $\Ha$-bimodule with
respect to the left and right coregular actions. We reserve for them the notation $\tr$ and $\tl$, respectively:
\be
x\tr a=a_{(1)}\:\langle a_{(2)},x\rangle, \quad
a\tl x=\langle a_{(1)},x\rangle\: a_{(2)},
\ee
where $x\in \Ha$, $a\in \Ha^*$, and $a_{(1)}\tp a_{(2)}=\Delta_{\Ha^*}(a)$.
$\Ha^*$ can be viewed as a left $\Ha'=\Ha^{op}\tp\Ha$-module algebra
with the universal R-matrix taken to be  $\Ru'=\Ru^{-1}_{13}\Ru_{24}$.
It is quasi-commutative, with condition (\ref{crel}) turning to the equation
\be{}
(a\tl{\Ru_1}) (b\tl{\Ru_2})  =   ({\Ru_2}\tr b) ({\Ru_1}\tr a)
\label{rtt}
\ee{}
for any $a,b\in \Ha^*$.
This equation is fulfilled due to (\ref{flip}).
\end{example}
\begin{example}
Along the line of Proposition \ref{main}, consider the twist from the Hopf algebra $\Ha^{op}\tp\Ha$ to
\tw{\Ha}{\Ru}{\Ha} by the cocycle $\Ru_{13}\Ru_{23}$.
Denote by $\tilde\Ha^*$ the  quasi-commutative \tw{\Ha}{\Ru}{\Ha}-module algebra
that is twist-equivalent to $\Ha^*$.
The equation
\begin{eqnarray}
({\Ru_{1'}}\tr a\tl{\Ru_2}) (b\tl{\Ru_1}\tl{\Ru_{2'}}) & = &
({\Ru_{1'}}\tr {\Ru_2}\tr b ) ({\Ru_1}\tr a\tl{\Ru_{2'}})
\label{re1}
\end{eqnarray}
taking place for any $a$ and $b$ from $\tilde\Ha^*$
is a consequence of (\ref{rtt}).
Here, the primes distinguish different copies of $\Ru$.
\end{example}
\end{parag}
\begin{parag}
The left action of  $\Ha^{op}\tp \Ha$ on $\Ha^*$
is expressed through the left and right coregular actions:
\be{}
\label{ttpaction}
(x\tp y)\btr a = y\tr a\tl \gm(x), \quad x\tp y \in\Ha^{op}\tp \Ha,\>a\in\Ha^*.
\ee{}
Let $m$ be the multiplication in  $\Ha^*$ and $\tilde m$ the multiplication
in $\tilde \Ha^*$.
They are related by the formulas
\begin{eqnarray}
\label{reprod}
\tilde m(a\tp b) =
 m\bigl(\Ru_1 \tr a \tl\Ru_{1'} \tp  b \tl \gm(\Ru_2) \tl\Ru_{2'} \bigr).\\
\label{frtprod}
  m(a\tp b) =
 \tilde m\bigl(\Ru_1 \tr a \tl\gm(\Ru_{1'}) \tp  b \tl \Ru_{2'} \tl\Ru_2 \bigr).
\end{eqnarray}
These expressions prove that the associative algebra $\tilde \Ha^*$ is exactly
the Hopf algebra in the quasitensor category of  $\Ha$-modules, \cite{Mj}.
\end{parag}
\begin{parag}[Universal K-matrix]
Let $\{e_i\}$ be a base in $\Ha$ and $\{e^i\}$ its dual in $\Ha^*$ with respect to
the canonical Hopf pairing. Let us consider
the element $\Ta=\sum_i e_i\tp e^i\in \Ha\tp \Ha^*$.
In case $\Ha$ is infinite dimensional this element can be defined by introducing an appropriate topology.
Using the identities
\be
\label{most_trivial}
\sum_i e_i\tp e^i \tl x = \sum_i x e_i\tp e^i
\quad
\sum_i e_i\tp x \tr e^i  = \sum_i e_i x\tp e^i,
\ee
which are true for any $x\in \Ha$, we rewrite
relations (\ref{rtt}) in a compact form that involves the universal T-matrix $\Ta$:
\be
\label{abs_rtt}
\Ru \Ta_1 \Ta_2 =  \Ta_2 \Ta_1 \Ru.
\ee
Now we consider the same base $\{e^i\}$ as that in $\tilde \Ha^*$. Then, the canonical element
$\K=\sum_i e_i\tp e^i \in \Ha\tp \tilde \Ha^*$ may be thought of
as a universal K-matrix. Indeed, applying identities (\ref{most_trivial})
to (\ref{re1}) we rewrite it in the form of the RE,
\be{}
\label{abs_re}
\Ru_{21} \K_1 \Ru \K_2 =  \K_2 \Ru_{21} \K_1 \Ru,
\ee{}
in $\Ha\tp \Ha \tp \tilde \Ha^*$.
\begin{remark}
In applications, elements of the dual base $\{e^i\}$ are usually expressed as polynomials
in generators of the algebra $\Ha^*$. Although belonging to the same linear space,
$e^i$ are expressed differently through generators of the deformed algebra $\tilde \Ha^*$.
While fixing the universal K-matrix as an element of $\Ha\tp \tilde \Ha^*$,
the proposed formula does not suggest an easy way of computing it.
\end{remark}
\end{parag}
\begin{parag}
\label{sUKM}
Let us consider an element $\K^\A \in \Ha \tp \A$ of the form
$\K^\A=\sum_i e_i\tp a^i$, were $\{e_i\}$ is a base in $\Ha$. It can be viewed as
a linear map $f\colon\tilde \Ha^*\to \A$, $f\colon u\mapsto \sum_i\langle u,e_i\rangle a^i.$
\begin{propn}
\label{fusion_prop}
The map $f\colon\tilde \Ha^*\to \A$ is an algebra homomorphism if and only if
\be{}
\label{fusion}
(\Delta\tp \id)\bigl(\K^\A\bigr)=\Ru^{-1}\K^\A_1\Ru \K^\A_2\in \Ha\tp\Ha\tp \A.
\ee{}
\end{propn}
\begin{proof}
Let us prove the statement for $\K^\A=\K$ first. This case corresponds to $\A=\tilde \Ha^*$ and
$f=\id$.
The universal T-matrix $\Ta\in \Ha\tp \Ha^*$ obeys the identity $(\Delta\tp\id)(\Ta)=\Ta_1\Ta_2$.
Now applying  formula (\ref{frtprod}) to the right-hand side of this equality we come to
equation (\ref{fusion}) for $\K^\A=\K$. The general case follows from here since
the element $\K^\A$ is equal to $(\id \tp f) (\K)$.
\end{proof}
\end{parag}
\begin{parag}[Characters of $\tilde \Ha^*$]
\label{sChar}
Let $\chi$ be a character (an algebra homomorphism to $\KK$) of the RE dual  $\tilde \Ha^*$.
We can think of it as an element from $\Ha$, then
$\chi=(\id \tp \chi)(\K)$; so $\rho(\chi)$ gives a
numeric solution to the RE in any representation $\rho$ of $\Ha$.
Using formula (\ref{fusion}) we find the necessary and
sufficient condition for an element $\chi \in \Ha$ to
be a character of $\tilde \Ha^*$:
\be
\label{char}
\Delta \chi= \Ru^{-1}\chi_1 \Ru \chi_2.
\ee
It is easy to see that this equation is preserved by
the similarity transformation $\chi\to \eta\chi \eta^{-1}$, where
$\eta$ is a group-like element or, in other words, a character
of the algebra $\Ha^*$ (the tensor square of a group-like element
commutes with $\Ru$).

The unit of $\Ha$  satisfies equation (\ref{char}), so the
counit $\ve_{\Ha^*}$ of $\Ha^*$ is a character of $\tilde \Ha^*$.
\end{parag}
\begin{parag}
The comultiplication map $\Delta \colon\Ha\to \Ha^{\otps 2}$
is a Hopf algebra embedding.
The composition
\be
(\Ru^{+}\tp \Ru^{-})\circ\Delta_{\Ha^*}\colon \Ha^{*op} \to \Ha^{\otps 2},
\ee
is a Hopf algebra homomorphism too. These
maps extend to a Hopf algebra homomorphism from the double $\mathrm{D}(\Ha)$
to $\Ha^{\otps 2}$, \cite{RS}. It follows that any $\Ha^{\otps 2}$-module
is that over $\Ha$, $\Ha^{*op}$, and  $\mathrm{D}(\Ha)$. This is true
for the algebra $\tilde \Ha^*$ in particular.

The element $\Q=\Ru_{21}\Ru\in \Ha\tp \Ha$ defines a map
\be
\label{Q}
\Q(a)=\langle a,\Q_1\rangle \Q_2, \quad a \in \tilde \Ha^*,
\ee
from $\tilde \Ha^*$ to $\Ha$. This map is an $\Ha$-equivariant
algebra homomorphism, where $\Ha$ is considered as the adjoint module
over itself\footnote{It appeared in \cite{FRT} as a map from
the FRT algebra, which is an extension of $\tilde \Ha^*$, to $\Ha$.}, \cite{Mj}.
The existence of this homomorphism follows from
the fact that the algebra $\tilde \Ha^*$ is
a left module algebra over $\Ha^{*op}$. The dual version of this statement claims that  $\tilde \Ha^*$ is
a right comodule algebra over $\Ha_{op}$. Then the map $\Q$ has the form
$(\ve_{\Ha^*}\tp \gm)\circ \delta$, where $\delta$ is the coaction
$\tilde\Ha^*\to \tilde\Ha^*\tp \Ha_{op}$, \cite{DM}:
\be{}
\label{Hcomod}
\delta(a)&=&a_{(2)}\tp
\langle a_{(3)},\Ru^{-1}_{1'}\rangle\langle \gm(a_{(1)}),\Ru_{2}\rangle
 \Ru^{-1}_{2'}\Ru_1,
\quad a\in \tilde \Ha^*.
\ee{}
Here the element $a$ is considered as belonging to the linear space
$\Ha^*\simeq \tilde \Ha^*$ and
$a_{(1)}\tp a_{(2)}\tp a_{(3)}$ is the symbolic decomposition of the
three-fold coproduct $\Delta^3_{\Ha^*}(a)$.
The counit $\ve_{\Ha^*}$ is a character of $\tilde\Ha^*$, cf. Subsection \ref{sChar}, and
this implies that $\Q$ is an algebra map. Clearly any character
of the associative algebra $\tilde \Ha^*$ defines in this way a
homomorphism from $\tilde \Ha^*$ to $\Ha$.

The element $\Q$ is the image of the universal K-matrix in $\Ha^{\tp2}$,
$\Q=(\id\tp \Q)(\K)$, so it satisfies the abstract reflection
equation in $\Ha^{\tp3}$, \cite{S}.
\end{parag}
\section{Braided  comultiplication and twist}
\label{sBCT}
\begin{parag}
\label{iota}
Let $\iota_k\colon\tilde \Ha^*\to \ttps{\tilde\Ha^*}{n}
$ be the embedding sending
$\tilde \Ha^*$ to the $k$-th tensor factor. It is an algebra homomorphism.
Consider the elements $\K^{\{k\}}=(\id \tp \iota_k)(\K)$.
\begin{propn}
The multiplication in $\ttps{\tilde\Ha^*}{n}$ satisfies
by the following commutation relations
\begin{equation}
\begin{array}{rcl}
\Ru_{21}\K^{\{k\}}_1 \Ru_{12} \K^{\{k\}}_2 &=&  \K^{\{k\}}_2 \Ru_{21} \K^{\{k\}}_1 \Ru_{12},\quad k=1,\ldots n,\\
\Ru^{-1}_{12} \K^{\{k\}}_1 \Ru_{12} \K^{\{l\}}_2  &=& \K^{\{l\}}_2 \Ru^{-1}_{12} \K^{\{k\}}_1\Ru_{12},\quad l,k=1,\ldots n,\>l<k.
\end{array}
\label{HH}
\end{equation}
\end{propn}
\begin{proof}
Denote by $\Ha^{\otps 2}_k$ the algebra $\Ha^{\otps 2}$ embedded in the $k$-th place in $(\Ha\otp \Ha)^{\otps n}$
and similarly $\tilde \Ha^*_k =\iota_k(\tilde \Ha^*)\subset \ttps{\tilde\Ha^*}{n}$.
The first line in (\ref{HH}) is the RE for $\K_k$, which holds since $\tilde \Ha^*_k$ is a subalgebra in
$\ttps{\tilde\Ha^*}{n}$.
The action of $\Ha^{\otps 2n}$ restricted to $\tilde \Ha^*_k$ is induced by that
of  $\Ha^{\otps 2}_k$ via the Hopf algebra projection $\Ha^{\otps 2 n}\to \Ha^{\otps 2}_k$.
Then the cross-relations between  $\tilde \Ha^*_k$ and  $\tilde \Ha^*_l$ are
 verified  straightforwardly.
\end{proof}
\begin{remark}
Relations (\ref{HH}) reflect the fact that $\ttps{\tilde\Ha^*}{n}$ is a quasi-commutative
$\Ha^{\otps2 n}$-module algebra, cf. Corollary \ref{qctp}.
\end{remark}
\end{parag}
\begin{parag}
In this subsection, we give an interpretation of the braided bialgebra structure
on $\tilde \Ha^*$ from  the Drinfeld's twist point of view.
Let us denote by $\Delta_{\tilde\Ha^*}$ the comultiplication  $\Delta_{\Ha^*}\colon \Ha^*\to \Ha^*\tp \Ha^*$
considered as a $\KK$-linear map $\tilde\Ha^*\to \mbox{\br{\tilde\Ha^*}{\Ru}{\tilde\Ha^*}}$.
\begin{thm}[\cite{Mj}]
\label{Mj}
The map  $\Delta_{\tilde\Ha^*}$ from $\tilde\Ha^*$ to the braided tensor product $\mbox{\br{\tilde\Ha^*}{\Ru}{\tilde\Ha^*}}$
 is
an algebra  homomorphism.
\end{thm}
\noindent
The proof of this theorem that is given in \cite{Mj} is quite complicated. Here we give another proof using Hopf algebra twist and
the universal K-matrix. The key observation is that
the braided tensor product of $\Ha^{\otps k}$- and $\Ha^{\otps m}$-module algebras, which is defined by (\ref{br_t_pr}),
is in fact a particular case of twisted tensor product from (\ref{br_pr}), as pointed out in Remark \ref{btp}, and hence  a module algebra over $\Ha^{\otps (k+m)}$.
Thus we reformulate Theorem \ref{Mj} in the following way.
\begin{thm}
For any $n\in \N$ the $n$-fold coproduct $\Delta^n_{\Ha^*}\colon \Ha^* \to \Ha^{*\tp n}$
considered as a $\KK$-linear map from $\tilde \Ha^*$ to the twisted tensor product $\ttps{\tilde\Ha^*}{n}$
is an algebra homomorphism.
\end{thm}
\begin{proof}
Put $\K^{\{1\ldots m\}}=\K^{\{1\}}\ldots\K^{\{m\}}\in \Ha\tp \ttps{\tilde\Ha^*}{m}$,  $m\in \N$.
A straightforward induction using formula (\ref{br_pr}) shows that
$\bigl(\id\tp\Delta^{m}_{\tilde\Ha^*}\bigr)(\K)= \K^{\{1\dots m\}}$ for all $m\in \N$.
We shall prove the theorem by induction on $n$.
The case $n=0$ holds because $\ve_{\Ha^*}$ is a character of $\tilde \Ha^*$, cf. Section \ref{sChar}, and the case $n=1$ is trivial.
Let us assume that the statement is true
for $n\geq 1$.
From the bottom line of equation (\ref{HH}), we find
\begin{equation}
\begin{array}{rcl}
\Ru^{-1} \K^{\{n+1\}}_1 \Ru \K^{\{1\ldots n\}}_2  &=& \K^{\{1\ldots n\}}_2 \Ru^{-1} \K^{\{n+1\}}_1\Ru.
\end{array}
\label{HHk}
\end{equation}
Put, in notation of Subsection \ref{sUKM}, $\A=\ttps{\tilde\Ha^*}{(n+1)}$
and $a^i=\Delta^{n+1}_{\tilde\Ha^*}(e^i)\in \A$.
Then, the element
$$
\K^\A=\K^{\{1\dots n+1\}}= \K^{\{1\dots n\}}\K^{\{n+1\}}
$$
fulfills condition (\ref{fusion}), as follows from the induction assumption and equation
(\ref{HHk}).
\end{proof}
\end{parag}
\begin{parag}
In this subsection we give, following \cite{KS}, a topological interpretation of the
algebra $\ttps{\tilde\Ha^*}{n}$.
It is known that a representation $(\rho,V)$ of $\Ha$ induces a representation of the braid group
$B_n$ of $n$ strands in $\R^3$.
Let $\{\si_i\}_{i=1}^{n-1}$ be the set of generators of $B_n$ satisfying
relations
$\si_{i} \si_{i+1}\si_{i}=\si_{i+1}\si_{i}\si_{i+1}$
for $i=1,\ldots,n-2$ and $\si_{i}\si_{j}=\si_{j}\si_{i}$ for $|j-i|>1$.
The representation of $B_n$ on $V^{\tp n}$ is defined via the correspondence
$\si_i\to S_{i,i+1}$,  where $S=PR$ is expressed through the image of the
 universal R-matrix $R=(\rho\tp\rho)(\Ru)$ and the flip operator $P$ on $V^{\tp 2}$.
As usually, the subscripts specialize embedding $\End(V)\subset \End^{\tp n}(V)$.

Let $A_m$ denote a set of $m$ parallel lines in $\R^3$.
The braid group $B^{m}_{n}$ of n strands in $\R^3\backslash A_m $
is generated by $\{\si_i\}_{i=1}^{n-1}\subset B_{n}\subset B_{n}^m$
and the elements $\{\tau_{k}\}_{k=1}^{m}$ subject
to relations
\be{}
\si_{n-1}\tau_k \si_{n-1}\tau_k&=&\tau_k \si_{n-1}\tau_k\si_{n-1},\quad k=1,\ldots,m,\\
\si^{-1}_{n-1}\tau_k \si_{n-1}\tau_l&=&\tau_l \si^{-1}_{n-1}\tau_k\si_{n-1},\quad k,l=1,\ldots,m,\>k>l.
\ee{}
Together with $\rho$, a representation $\theta$ of the algebra $\ttps{\tilde\Ha^*}{m}$ on
a module $W$ defines a representation of $B_n^m$ on the module $V^{\tp n}\tp W$.
It extends the representation of the subgroup $B_n\subset B_n^m$ in $V^{\tp n}$ via
the natural embedding $V^{\tp n}\subset V^{\tp n}\tp W$.
Further, assigned to every $\tau_k$ is the operator $(\rho\tp\theta)\circ (\id \tp\iota_k)(\K)$,
where the embedding $\iota_k$ is defined in Subsection \ref{iota}.
\end{parag}
\section{Fusion rules for RE matrices}
\label{sFusion}
\begin{parag}
\label{fusion_rules}
In this section we consider the following problem.
Let $\bigl(\rho^{\{1\}},V^{\{1\}}\bigr)$ and $\bigl(\rho^{\{2\}},V^{\{2\}}\bigr)$ be two representations of $\Ha$.
Let $\A$ be an associative algebra and $K^{\{i\}}$, $i=1,2$, two constant RE
matrices with coefficients in $\A$ in the representations $\rho^{\{1\}}$ and $\rho^{\{2\}}$.
How to build  an RE matrix in the tensor product representation $\rho^{\{1\}}\tp \rho^{\{2\}}$ or
in its sub-representation?
\end{parag}
\begin{parag}
Let $\rho^{\{i\}}$ be representations of $\Ha$
on modules $V^{\{i\}}$, where $i$ runs over a set of indices, $i\in\I$.
Denote by $R^{\{i,j\}}=\bigl(\rho^{\{i\}}\tp\rho^{\{j\}}\bigr)(\Ru)$
the image of the universal R-matrix in $\End\bigl(V^{\{i\}}\bigr)\tp\End\bigl(V^{\{j\}}\bigr)$,
$i,j\in\I$.
These matrices satisfy the Yang-Baxter equation
$$
R^{\{i,j\}}_{\;\;12}R^{\{i,k\}}_{\;\;13}R^{\{j,k\}}_{\;\;23}=
R^{\{j,k\}}_{\;\;23}R^{\{i,k\}}_{\;\;13}R^{\{i,j\}}_{\;\;12}, \quad i,j,k\in\I.
$$
Let $K^{\{i\}}\in \End\bigl(V^{\{i\}}\bigr)\tp \A$ be  RE matrices in representations $\rho^{\{i\}}$, $i\in \I$,
and consider the elements
\be
\label{fusion_mat}
K^{\{i,j\}}=\bigl(R^{\{i,j\}}\bigr)^{-1}K^{\{i\}}_1R^{\{i,j\}}K^{\{j\}}_2\in \End\bigl(V^{\{i\}}\tp V^{\{j\}}\bigr)\tp \A,
\quad i,j\in \I.
\ee
\begin{propn}
\label{fusion_prop1}
Suppose the RE matrices have the form
$K^{\{i\}}=\bigl(\rho^{\{i\}}\tp f\bigr)(\K)$, $i\in \I$, where $f$ is an algebra homomorphism $\tilde \Ha^*\to \A$.
Then, the elements $K^{\{i,j\}}$ defined by (\ref{fusion_mat})
are RE matrices in the tensor product representation $\rho^{\{i\}}\tp \rho^{\{j\}}$.
They can be restricted to any sub-representation.
\end{propn}
\begin{proof}
Apply the algebra homomorphism $\bigl(\rho^{\{i\}}\tp \rho^{\{j\}}\tp f\bigr)\circ (\Delta\tp \id)$ to the
universal K-matrix and use Proposition \ref{fusion_prop}.
\end{proof}
\end{parag}
\begin{parag}
Proposition \ref{fusion_prop1} allows to build tensor product of  solutions to the RE
of special type, namely when they are images of the same universal K-matrix (cf. concluding remarks
below) via the
same homomorphism $f$. Since {\em a priori} there are no  cirteria for this
to be true\footnote{Even in the case $\A=\C$ and the fundamental vector representation
of $\U_q\bigl(gl(n)\bigr)$ there is a variety of solutions to the same
RE, \cite{KSS,M}.}, one has to substitute this condition by certain compatibility requirements.
\begin{definition}
An {\em RE data} is a set $\bigl\{\bigl(\rho^{\{i\}},V^{\{i\}},K^{\{i\}}\bigr)\bigr\}_{i\in \I}$
of triples, where $\rho^{\{i\}}$ are finite dimensional representations of $\Ha$ on $V^{\{i\}}$ and
the matrices $K^{\{i\}}\in \End\bigl(V^{\{i\}}\bigr)\tp \A$  satisfy the equations
\be
\label{compatibility}
R^{\{j,i\}}_{\;\;21}K^{\{i\}}_1R^{\{i,j\}}K^{\{j\}}_2=K^{\{j\}}_2R^{\{j,i\}}_{\;\;21}K^{\{i\}}_1R^{\{i,j\}}
\ee
in $\End\bigl(V^{\{i\}}\tp V^{\{j\}}\bigr)\tp \A$, $i,j\in \I$.
\end{definition}
\noindent
Clearly restriction of the set $\I$ to any its subset gives an RE data.
In particular, each triple  $\bigl(\rho^{\{i\}},V^{\{i\}},K^{\{i\}}\bigr)$
is an RE data for $\I=\{i\}$; condition (\ref{compatibility})
reduces to the RE equation on $K^{\{i\}}$.
Given an RE triple $(\rho,V,K)$ and a set $\I$ one can form an RE data
putting $\bigl(\rho^{\{i\}},V^{\{i\}},K^{\{i\}}\bigr)=(\rho,V,K)$ for each $i\in \I$.
Also, for any homomorphism $f\colon \tilde \Ha^*\to \A$
the matrices $K^{\{i\}}=\bigl(\rho ^{\{i\}}\tp f\bigr)(\K)$ satisfy condition (\ref{compatibility}).
\begin{propn}
\label{restriction1}
Consider a triple $(\rho,V,K)$ from an RE data $\mathfrak{K}$. Suppose the
module $V$ is semisimple. Let $(\rho_0,V_0)$ be a sub-representation of
$(\rho,V)$ and $K_0$ the restriction of $K$ to $V_0$ from Proposition \ref{restriction}.
Then, the set $\mathfrak{K}\cup \left\{(\rho_0,V_0,K_0)\right\}$ is an RE data.
\end{propn}
\begin{proof}
The apparent modification of the proof of Proposition \ref{restriction}.
\end{proof}
\noindent
The following theorem generalizes Proposition \ref{fusion_prop1}.
\begin{thm}[Fusion procedure]
\label{fusion_thm}
Let $\bigl\{\bigl(\rho^{\{i\}},V^{\{i\}},K^{\{i\}}\bigr)\bigr\}_{i\in \I}$ be an RE data.
Then, the union
\be
\label{union}
\bigl\{\bigl(\rho^{\{i\}},V^{\{i\}},K^{\{i\}}\bigr)\bigr\}_{i\in \I} \cup
\;\bigl\{\bigl(\rho^{\{i\}}\tp\rho^{\{j\}},V^{\{i\}}\tp V^{\{j\}},K^{\{i,j\}}\bigr)\bigr\}_{(i,j)\in \I\times \I}\;,
\ee
where $K^{\{i,j\}}$ is defined in (\ref{fusion_mat}), is an RE data, too.
\end{thm}
\begin{proof}
See Appendix.
\end{proof}
\begin{remark}
It is easy to prove "associativity" of the fusion: $K^{\{\{i,j\},k\}}=K^{\{i,\{j,k\}\}}$.
\end{remark}
\noindent
Let us apply Proposition \ref{restriction1} and Theorem \ref{fusion_thm} to the quantum group $\U_q(\g)$, where $\g$ is a semisimple
Lie algebra. Its finite dimensional modules are semisimple, and each irreducible representation
can be realized as  a submodule in the tensor power of a fundamental representation, \cite{ChPr}.
Starting from an RE matrix in any fundamental representation,
one can build an RE matrix in every irreducible representation by tensoring and projecting.
\end{parag}
\section{Concluding remarks}
Suppose a quasitriangular Hopf algebra $\Ha$ is a subbialgebra
in a bialgebra $\B$ which is quasitriangular
in the sense that condition (\ref{flip}) is fulfilled for any
$x\in \B$. One can take as $\B$ the dual to
the FRT bialgebra associated with a finite dimensional representation
of $\Ha$. All the constructions of this paper can be literally carried
over  from the dual Hopf algebra $\Ha^*$ to the bialgebra
$\B^*$, which is  a quasi-commutative
bimodule algebra over $\Ha$. One can twist $\B^*$ along the line of
Proposition \ref{main} to obtain an \tw{\Ha}{\Ru}{\Ha}-module
algebra, $\tilde \B^*$.
There is a braided bialgebra structure on $\tilde\B^*$ given
by the coproduct of $\B^*$, which is considered as a map from
$\tilde\B^*$ to $\ttps{\tilde\B^*}{2}$.
The algebra $\tilde \B^*$ is a comodule over
$\Ha^{*op}$ and the counit $\ve_{\B^*}$ is its character.
So there is an $\Ha$-equivariant homomorphism
$\Q\colon \tilde \B^*\to \Ha$ defined by  formula
(\ref{Q}) where $a \in \tilde\B^*$.
Likewise for $\B=\Ha$, one can introduce
the universal K-matrix $\K\in \B\tp \tilde \B^*$, which is
the canonical element of the bialgebra pairing
between $\B$ and $\B^*$. It will satisfy
equations (\ref{abs_re}) and (\ref{fusion}) in $\B\tp\B\tp\tilde\B^*$.
Any representation of $\B$ gives a solution to the corresponding
RE matrix, and the formula (\ref{fusion_mat}) provides
a fusion procedure for such RE matrices.

There is an algebra map $\tilde \B^*\to \tilde\Ha^*$ coming
from the bialgebra homomorphism $\B^*\to \Ha^*$; the latter is
dual to the embedding $\Ha\to \B$.
Therefore $\tilde \B^*$ gives, in general, more solutions to the RE than $\tilde \Ha^*$
does. If, for instance, $\Ha$ is a factorizable Hopf algebra, \cite{RS},
the map $\Q\colon \tilde\Ha^*\to \Ha$ is an algebra isomorphism.
So the characters of $\tilde \Ha^*$ are exactly those of $\Ha$.
Consider, e.g., the quantum group $\U_q\bigl(gl(n)\bigr)$.
It has only one-parameter family of one-dimensional representations
assigning a scalar to the central generator.
On the other hand, there are much more solutions
to the RE in the fundamental vector representation of  $\U_q\bigl(gl(n)\bigr)$
with the same R-matrix, \cite{KSS,M}. They are characters
of $\tilde\B^*$, which is in this case the matrix RE algebra
associated with the given representation of $\U_q\bigl(gl(n)\bigr)$.

\section{Appendix: proof of Theorem \ref{fusion_thm}}
\begin{parag}
First we claim that for every $i,j,k\in \I$ the following identities hold true:
\be
R^{\{k,j\}}_{\;\;32}R^{\{k,i\}}_{\;\;31}
\left\{\bigl(R^{\{i,j\}}_{\;\;12}\bigr)^{-1} K^{\{i\}}_1R^{\{i,j\}}_{\;\;12}K^{\{j\}}_2\right\}
R^{\{i,k\}}_{\;\;13}R^{\{j,k\}}_{\;\;23}
\left\{ K^{\{k\}}_3\right\}=
\nonumber\hspace{3cm}\\\hspace{3cm}
=\left\{ K^{\{k\}}_3\right\}R^{\{k,j\}}_{\;\;32}R^{\{k,i\}}_{\;\;31}
\left\{\bigl(R^{\{i,j\}}_{\;\;12}\bigr)^{-1} K^{\{i\}}_1R^{\{i,j\}}_{\;\;12}K^{\{j\}}_2\right\}
R^{\{i,k\}}_{\;\;13}R^{\{j,k\}}_{\;\;23},
\label{partial1}
\\
R^{\{j,i\}}_{\;\;21}R^{\{k,i\}}_{\;\;31}\left\{K^{\{i\}}_1\right\}
R^{\{i,k\}}_{\;\;13}R^{\{i,j\}}_{\;\;12}
\left\{\bigl(R^{\{j,k\}}_{\;\;23}\bigr)^{-1} K^{\{j\}}_2R^{\{j,k\}}_{\;\;23}K^{\{k\}}_3\right\}=
\nonumber\hspace{3cm}\\\hspace{3cm}
\left\{\bigl(R^{\{j,k\}}_{\;\;23}\bigr)^{-1} K^{\{j\}}_2R^{\{j,k\}}_{\;\;23}K^{\{k\}}_3\right\}
R^{\{j,i\}}_{\;\;21}R^{\{k,i\}}_{\;\;31}\left\{K^{\{i\}}_1\right\}
R^{\{i,k\}}_{\;\;13}R^{\{i,j\}}_{\;\;12}.
\label{partial3}
\ee
Let us prove equation (\ref{partial1}).
Pulling $R^{\{k,j\}}_{\;\;32}$ from left to right in the l.h.s. of (\ref{partial1}) and using the Yang-Baxter equation
twice we obtain
$$
\bigl(R^{\{i,j\}}_{\;\;12}\bigr)^{-1}R^{\{k,i\}}_{\;\;31}K^{\{i\}}_1
R^{\{i,k\}}_{\;\;13}
R^{\{i,j\}}_{\;\;12}
R^{\{k,j\}}_{\;\;32}
K^{\{j\}}_2
R^{\{j,k\}}_{\;\;23}
K^{\{k\}}_3.
$$
Now pulling $K^{\{k\}}_3$ from right to left and employing (\ref{compatibility}) twice,
we get
$$
K^{\{k\}}_3
\bigl(R^{\{i,j\}}_{\;\;12}\bigr)^{-1}R^{\{k,i\}}_{\;\;31}K^{\{i\}}_1
R^{\{i,k\}}_{\;\;13}
R^{\{i,j\}}_{\;\;12}
R^{\{k,j\}}_{\;\;32}
K^{\{j\}}_2
R^{\{j,k\}}_{\;\;23}.
$$
This time we pull $R^{\{k,j\}}_{\;\;32}$ from right to to left and apply
the Yang-Baxter equation twice to obtain the r.h.s. of (\ref{partial1}).

Now we turn to equation  (\ref{partial3}). Pulling $\bigl(R^{\{j,k\}}_{\;\;23}\bigr)^{-1}$ to the left in the l.h.s. and
using the Yang-Baxter equation twice we come to
$$
\left\{\bigl(R^{\{j,k\}}_{\;\;23}\bigr)^{-1}\right\}R^{\{k,i\}}_{\;\;31}R^{\{j,i\}}_{\;\;21}\left\{K^{\{i\}}_1\right\}
R^{\{i,j\}}_{\;\;12}R^{\{i,k\}}_{\;\;13}
\left\{ K^{\{j\}}_2R^{\{j,k\}}_{\;\;23}K^{\{k\}}_3\right\}.
$$
This is equal
$$
\left\{\bigl(R^{\{j,k\}}_{\;\;23}\bigr)^{-1} K^{\{j\}}_2\right\}R^{\{k,j\}}_{\;\;31}R^{\{j,i\}}_{\;\;21}\left\{K^{\{i\}}_1\right\}
R^{\{i,j\}}_{\;\;12}R^{\{i,k\}}_{\;\;13}
\left\{R^{\{j,k\}}_{\;\;23}K^{\{k\}}_3\right\},
$$
due to (\ref{compatibility}). Pulling $R^{\{j,k\}}_{\;\;23}$  to the left and apply the Yang-Baxter equation twice
we obtain
$$
\left\{\bigl(R^{\{j,k\}}_{\;\;23}\bigr)^{-1} K^{\{j\}}_2R^{\{j,k\}}_{\;\;23}\right\}
R^{\{j,i\}}_{\;\;21}R^{\{k,i\}}_{\;\;31}\left\{K^{\{i\}}_1\right\}
R^{\{i,k\}}_{\;\;13}R^{\{i,j\}}_{\;\;12}
\left\{K^{\{k\}}_3\right\}.
$$
Now we push $K^{\{k\}}_3$  to the left, employ (\ref{compatibility}), and come to the r.h.s. of (\ref{partial3}).
\end{parag}
\begin{parag}
Now we prove that the matrix $K^{\{i,j\}}$, $i,j\in \I$, obeys the reflecion equation in
the representation $\rho^{\{i\}}\tp \rho^{\{j\}}$.
The image of the universal R-matrix in the representation $\rho^{\{i\}}\tp \rho^{\{j\}}$
is equal to
$$R^{\{i,j\}}_{\;\;14}R^{\{i,i\}}_{\;\;13}R^{\{j,j\}}_{\;\;24}R^{\{j,i\}}_{\;\;23}\in
\left\{\End\bigl(V^{\{i\}}\bigr)\tp\End\bigl(V^{\{j\}}\bigr)\right\}\tp \left\{\End\bigl(V^{\{i\}}\bigr)\tp\End\bigl(V^{\{j\}}\bigr)\right\}.$$
The left-hand side of the RE on $K^{\{i,j\}}$
reads
\be
R^{\{i,j\}}_{\;\;32}R^{\{i,i\}}_{\;\;31}R^{\{j,j\}}_{\;\;42}R^{\{j,i\}}_{\;\;41}
\left\{\bigl(R^{\{i,j\}}_{\;\;12}\bigr)^{-1} K^{\{i\}}_1R^{\{i,j\}}_{\;\;12}K^{\{j\}}_2\right\}
\times
\nonumber\hspace{3cm}\\\hspace{3cm}
\times
R^{\{i,j\}}_{\;\;14}R^{\{i,i\}}_{\;\;13}R^{\{j,j\}}_{\;\;24}R^{\{j,i\}}_{\;\;23}
\left\{\bigl(R^{\{i,j\}}_{\;\;34}\bigr)^{-1} K^{\{i\}}_3R^{\{i,j\}}_{\;\;34}K^{\{j\}}_4\right\},
\nonumber
\ee
while the r.h.s. of the RE is obtained from this by placing the last factor concluded in $\{\;\}$ to  the leftmost
position. Pulling $\bigl(R^{\{i,j\}}_{\;\;34}\bigr)^{-1}$ to the left and using
the Yang-Baxter equation four times, we obtain
\be
\left\{\bigl(R^{\{i,j\}}_{\;\;34}\bigr)^{-1}\right\}
R^{\{j,j\}}_{\;\;42}R^{\{j,i\}}_{\;\;41}R^{\{i,j\}}_{\;\;32}R^{\{i,i\}}_{\;\;31}
\left\{\bigl(R^{\{i,j\}}_{\;\;12}\bigr)^{-1} K^{\{i\}}_1R^{\{i,j\}}_{\;\;12}K^{\{j\}}_2\right\}
\times
\nonumber\hspace{3cm}\\\hspace{3cm}
\times
R^{\{i,i\}}_{\;\;13}R^{\{j,i\}}_{\;\;23}R^{\{i,j\}}_{\;\;14}R^{\{j,j\}}_{\;\;24}
\left\{ K^{\{i\}}_3R^{\{i,j\}}_{\;\;34}K^{\{j\}}_4\right\}.
\nonumber
\ee
From this, using (\ref{partial1}) for $k=i$, we come to
\be
\left\{\bigl(R^{\{i,j\}}_{\;\;34}\bigr)^{-1}K^{\{i\}}_3\right\}
R^{\{j,j\}}_{\;\;42}R^{\{i,j\}}_{\;\;32}R^{\{j,i\}}_{\;\;41}R^{\{i,i\}}_{\;\;31}
\left\{\bigl(R^{\{i,j\}}_{\;\;12}\bigr)^{-1} K^{\{i\}}_1R^{\{i,j\}}_{\;\;12}K^{\{j\}}_2\right\}
\times
\nonumber\hspace{3cm}\\\hspace{3cm}
\times
R^{\{i,i\}}_{\;\;13}R^{\{i,j\}}_{\;\;14}R^{\{j,i\}}_{\;\;23}R^{\{j,j\}}_{\;\;24}
\left\{ R^{\{i,j\}}_{\;\;34}K^{\{j\}}_4\right\}.
\nonumber
\ee
Now we push
$R^{\{i,j\}}_{\;\;34}$ to the left and use the Yang-Baxter equation four times.
The resulting expression is
\be
\left\{\bigl(R^{\{i,j\}}_{\;\;34}\bigr)^{-1}K^{\{i\}}_3 R^{\{i,j\}}_{\;\;34}\right\}
R^{\{i,j\}}_{\;\;32}R^{\{i,i\}}_{\;\;31}R^{\{j,j\}}_{\;\;42}R^{\{j,i\}}_{\;\;41}
\left\{\bigl(R^{\{i,j\}}_{\;\;12}\bigr)^{-1} K^{\{i\}}_1R^{\{i,j\}}_{\;\;12}K^{\{j\}}_2\right\}
\times
\nonumber\hspace{1.5cm}\\\hspace{3cm}
\times
R^{\{i,j\}}_{\;\;14}R^{\{j,j\}}_{\;\;24}R^{\{i,i\}}_{\;\;13}R^{\{j,i\}}_{\;\;23}
\left\{ K^{\{j\}}_4\right\}.
\nonumber
\ee
To accomplish the proof, one should push $K^{\{j\}}_4$ to the left employing the identity
(\ref{partial1}) where $k$ is equal to $j$ and $3$ is replaced by $4$.
\end{parag}
\begin{parag}
It remains to check that for any choice of $i,j,k\in \I$, the matrices $K^{\{i,j\}}$, $K^{\{k\}}$
are compatible  in the sense of (\ref{compatibility}). This statement splits into several assertions.
First of all, the matrices $K^{\{i,j\}}$ and $K^{\{k\}}$ must solve the RE in their own
representations. In what concerns  $K^{\{k\}}$, this holds by assumption;
as to the matrix $K^{\{i,j\}}$, this has been proven in the previous subsection.
We must verify the cross-relations between  $K^{\{i,j\}}$ and $K^{\{k\}}$.
It is not difficult to see that they are encoded in  the system of equations (\ref{partial1}--\ref{partial3}).
\end{parag}

\bigskip

\end{document}